\documentstyle[11pt,twoside, amscd]{amsart}

\newcommand{\cT}{{\cal T}}

\newcommand{\cC}{{\cal C}}

\newcommand{\cM}{{\cal M}}

\newcommand{\cA}{{\cal A}}
\newcommand{\cF}{{\cal F}}

\newcommand{\Spec}{\operatorname{Spec}}

\newcommand{\Mor}{\operatorname{Mor}}
\newcommand{\Ch}{\operatorname{Chow}}

\newcommand{\gene}{\operatorname{span}}

\newcommand{\relint}{\operatorname{relint}}

	\newcommand{\Hom}{\operatorname{Hom}}

        \newfont{\hollow}{msbm10 scaled\magstep1}
        \newfont{\Bfmit}{eufm10 scaled\magstep1}

\def\proof{\smallskip\noindent{\it Proof. }}
\def\endproof{\hfill\qed}
\def\qed {\nobreak$\quad$\lower 1pt\vbox{
    \hrule
    \hbox to 8pt{\vrule height 8pt\hfil\vrule height 8pt}
      \hrule}\ifmmode\relax\else\par\medbreak \fi}

\newtheorem{thm}{Theorem}[section]

\newtheorem{question}[thm]{Question}
\newtheorem{prob}[thm]{Problem}
\newtheorem{lem}[thm]{Lemma}
\newtheorem{cor}[thm]{Corollary}

\newtheorem{prop}[thm]{Proposition}

\newtheorem{conj}[thm]{Conjecture}
\theoremstyle{definition}
\newtheorem{defn}[thm]{Definition}

\newtheorem{exmp}[thm]{Example}

\newtheorem{rem}[thm]{Remark}           

\theoremstyle{remark}

\begin{document}
\title[Combinatorics and Quotients of Toric Varieties]{Combinatorics and  Quotients of Toric Varieties}
\author[Yi {\sc Hu}]{Yi {\sc Hu}}

\maketitle

\section{Introduction}
This paper studies two related subjects. One is some combinatorics arising
from linear projections of polytopes and fans of cones. The other
is quotient varieties of toric varieties. The relation is that projections of
polytopes are related to  quotients of projective toric varieties and
projection of fans are related to quotients of general toric varieties.
Despite its relation to geometry
the first part is  purely combinatorial and should  be  of interest
in  its own right.

For the combinatorial part, consider
a linear projection $\pi: V \rightarrow W$
between two real vector spaces.  Let $P$ be a polytope in $V$ of full dimension and $Q = \pi (P)$.
The projections of all faces of $P$ under $\pi$ 
induce a polytopal subdivision of $Q$. This gives rise to 
the poset $\Gamma = \Gamma (P,Q)$ of the so-called {\it chambers} in $Q$ (Definition 1.1).

The sequence $0 \rightarrow \ker \pi \rightarrow V @>\pi>> W \rightarrow 0$
induces the dual sequence
$$0 \rightarrow W^* \rightarrow V^* @>\pi^\vee>> (\ker \pi)^* \rightarrow 0.$$
The projections of all cones of the normal fan $\Delta (P)$ of $P$ under $\pi^\vee$
induce a fan in $(\ker \pi)^*$. 
This is  the normal fan $\Delta (\Sigma (P,Q))$
of the Billera-Sturmfels fiber polytope $\Sigma (P,Q)$. Let
 $\Gamma^* = \Gamma^*(P,Q)$ be the poset of cones in this fan (Definition 1.2).

The above descriptions clearly put $\Gamma$ and $\Gamma^*$ in a linear dual situation. Indeed,
they are part of a larger picture which goes beyond just dual analogy.

First,  it has been known that $\Gamma^*$
is in one-to-one correspondence with the poset $\cT_{coh} = \cT_{coh} (P, Q)$ of the so-called
{\it coherent strings} of the projection $\pi: P \to Q$ (Definition 2.1). 
The poset $\cT = \cT (P, Q)$ of {\it locally coherent strings},
 the generalizations of coherent strings, has also been around for quite a while (Definition 3.1). Investigating
into this, it is quite natural to ask the following questions:
\begin{enumerate}
\item what are the generalizations of elements in $\Gamma^*$ so that  they correspond to locally coherent        strings in $\cT$?
\item  are there {\it coherent costrings} linear dual to coherent strings?
\item  are there {\it locally coherent costrings} linear dual to  locally coherent strings?
\item  are there perfect correspondences for these new objects similar to the correspondence between
      $\Gamma^*$ and $\cT_{coh}$?
\end{enumerate}

Indeed, perfect correspondences do occur and are the main theme of this paper.
 To this end, some new objects are to be introduced
in response to each and every of the above questions:
\begin{enumerate}
\item the poset $\Gamma^*_{vir} = \Gamma^*_{vir} (P, Q)$ of {\it virtual cones} of  the projection $\pi^\vee$.
     A virtual cone is a combinatorial abstraction that generalizes a real polyhedral cone
  (Definition 7.1). 
     It corresponds to a locally coherent string of $\pi: P \to Q$.
\item the poset $\cT^*_{coh}= \cT^*_{coh} (P, Q)$ 
      of {\it coherent costrings} of the projection $\pi^\vee$. This notion is linear
     dual to that of coherent strings. A coherent costring is roughly speaking a {\it special}
     lifting of a fan via the map $\pi^\vee$ (Definition 4.2). It corresponds to a chamber in $\Gamma$.
\item the poset of $\cT^* = \cT^* (P, Q)$ of {\it locally coherent costrings} of the projection $\pi^\vee$. 
This notion is linear dual to that of locally coherent strings. A locally
coherent costring is roughly speaking a lifting of a fan via the map $\pi^\vee$ (Definition 5.1).
\item the poset $\Gamma_{vir} = \Gamma_{vir} (P, Q)$ of {\it virtual chambers (cells)} of the projection $\pi$.
 A virtual chamber (cell) is a combinatorial abstraction that 
generalizes a real chamber (cell) in $\Gamma$ (Definition 6.1).
 It corresponds to a locally coherent costring.
\end{enumerate}

Locally coherent strings first appeared in \cite{BKS94} when $\dim Q=1$ and were generalized
and studied afterwards, especially in association with the {\it generalized Baues conjecture}
 (see \cite{Reiner98} for a recent survey).
Virtual chambers for  a simplex was introduced and studied in \cite{dHSS}.

Our first main result in this paper is
\proclaim Theorem A. 
\label{main} (Theorem 8.1.) 
 Let the notations be as explained in the above. Then 
\begin{enumerate}
\item  $\cT (P,Q)$ is canonically anti-isomorphic to $\Gamma^*_{vir} (P, Q)$ which extends the 
anti-isomorphism between $\cT_{coh}(P, Q) $ and $\Gamma^* (P, Q)$;
\item  $\cT^* (P, Q)$ is canonically anti-isomorphic to $\Gamma_{vir} (P, Q)$
which extends the anti-isomorphism between $\cT^*_{coh} (P, Q) $ and $\Gamma (P, Q)$. 
\end{enumerate}

As hinted in the beginning, the above combinatorial objects underlie some
quotient toric geometry. We now elaborate. 

Let $X_{\Delta_0} = \cup_{\sigma \in \Delta_0} A_\sigma$ be the toric variety over the complex
number field
defined by a fan $\Delta_0$ with a lattice $N$,
 where  $A_\sigma$ is the affine variety defined by the cone $\sigma$ in $\Delta_0$.  
Let $M^2$ be any sublattice of the dual lattice $M$ of $N$. Consider the exact sequence
$$0 \rightarrow  M^2 \rightarrow  M @>{\pi}>>  M^1 \rightarrow  0.$$
If we apply $\Hom_{\Bbb Z} (-, {\Bbb C}^*)$ to it, then we get the exact sequence
$$1  \rightarrow \Hom_{\Bbb Z} (M^1, {\Bbb C}^*) \rightarrow \Hom_{\Bbb Z} (M, {\Bbb C}^*) \rightarrow \Hom_{\Bbb Z} (M^2, {\Bbb C}^*) \rightarrow 1.$$
The group $G = \Hom_{\Bbb Z} (M^1, {\Bbb C}^*)$ is the product of the torus 
$\Hom_{\Bbb Z} (M^1_{free}, {\Bbb C}^*)$ and the finite group
$\Hom_{\Bbb Z} (M^1_{tor}, {\Bbb C}^*)$.
Let $$0 \rightarrow  N^1 \rightarrow  N @>{\pi^\vee}>> N^2 \rightarrow  0$$
be the dual sequence.

The {\it projective} quotients of {\it projective} toric varieties by subtorus
have been studied by Kapranov, Sturmfels and Zelevinsky in \cite{KSZ}. See also, \cite{Hu91}.
For the arbitrary toric variety $X_{\Delta_0}$ (not necessarily complete, not necessarily quasi-projective),
we have the following 

\proclaim Theorem B.
 (Theorem 11.5)
Let $\Delta \subset \Delta_0$ be a subset of cones and $U = \cup_{\sigma \in \Delta} A_\sigma$ be an open
subset of $X_{\Delta_0}$ where $A_\sigma$ is the affine open subset defined by $\sigma$. Then
\begin{enumerate}
\item $U$ has a categorical quotient by the action of $G$ if $\Delta$ is a locally coherent costring of the projection $\pi^\vee$
and every cone in $\pi^\vee (\Delta)$ is strongly convex. 
    In this case the quotient
   $U /\!/G $ is isomorphic to $X_{\pi^\vee (\Delta)}$, the toric variety defined by the induced fan
   $\pi^\vee (\Delta)$;
\item $U /\!/G$ is geometric if and only if $\dim \pi^\vee (\sigma) = \dim \sigma $ for all $\sigma \in \Delta$
(i.e., $\Delta$ is tight by definition, or equivalently, minimal in $\cT^*$).
\end{enumerate}

When $X_{\Delta_0}$ is projective and equipped with a $T$-linearized ample line bundle $L$,
there corresponds to a polytope $P$ in the  lattice $M$. Let $Q= \pi (P)$.
Then by the equivalence of (2)  in Theorem A, we have

\proclaim Corollary C. Every virtual chamber (cell) in $\Gamma_{vir} (P,Q)$
defines a (not-necessarily projective) toric quotient variety of $X_{\Delta_0}$
whose fan is induced from the corresponding locally coherent 
costring in $\cT^* (P, Q)$.

 Note that the real chambers (cells) correspond to projective toric quotients
whose fans are those given by global coherent costrings.

In the end of this paper, we try to relate the combinatorial bistellar flips to the
geometric Mori-type flips among quotients of affine spaces.

Also we introduce a new notion, {\it virtual} (or {\it pseudo}) oriented matroid, 
which is somewhat related to the rest of the paper.


Finally, we mention that there are some recent works that study the quotients of toric varieties
by subtori  \cite{BBSw1}, \cite{BBSw2}, \cite{ac1} and \cite{ac2}.
Note that the group $G$ that we use to take quotient  is typically disconnected. 
Also, I learned first from Reiner  that  Santos (\cite{Santos})
has obtained some important  results on Baues conjectures which
are related to our combinatorial constructions. 

{\bf Acknowledgments.} 
My thanks go to Bernd Sturmfels for suggesting that I investigate the perfect correspondences
that appeared in this paper.  This paper grew out of our conversation at Utah Algebraic
Geometry Conference in 1995 and thereafter and I thank him for his  suggestions, comments and
encouragements.  I owe thanks to Victor Reiner for saving me from an embarrassment.
I thank Michel Brion  and G\"unter Ziegler for some useful conversations or correspondence
regarding the materials in this paper, and also Andrzej Bialynicki-Birula for sending his papers
to me.

{\scriptsize
\tableofcontents
}

\section{Chambers and cones}

Let $ V$ and $W$ be vector spaces over
${\Bbb R}$ of finite dimensions, and $$\pi: V \rightarrow W$$ be a linear projection.
Assume that $P \subset V$ is a polytope of full dimension and $Q = \pi(P)$.

\begin{defn}
The (polytopal) cell in $Q$ containing a point $q \in Q$ is the intersection of 
$\pi$-images of all the faces of $P$
that contain the given point $q$, that is
$${\bf c} (q) :=[q] :=  \cap \{\pi(F) | q \in \pi(F), F \in L(P) \}.$$
The set of all cells is denoted by $\Gamma = \Gamma (P, Q)$,
 partially ordered by inclusion. Maximal cells are called
 {\it chambers}\footnote{Cells are often called chambers in other literature. Then,
chambers in our sense would have to be called maximal chambers.}.
\end{defn}

Dually, the sequence $0 \rightarrow \ker \pi \rightarrow V @>{\pi}>>  W \rightarrow 0$
induces the sequence
$$0 \rightarrow W^* \rightarrow V^* @>{\pi^\vee}>> (\ker \pi)^* \rightarrow 0.$$

 Recall that given any linear function $\psi \in V^*$, the set $P^\psi$ of all points in $P$
on which $\psi$ is maximal is a face of $P$, and all nonempty faces of $P$ arise in this way.
In fact, given any $F$ in the {\it face lattice} $L(P)$ of $P$, the set of $\psi \in V^*$ such that
$P^\psi = F$ is a cone, called  the normal cone of $F$, and denoted by  $N(P, F)= N(F)$.
The set  $\Delta (P)$ of the normal cones of all faces of $P$ is the normal fan of $P$.

\begin{defn}
\label{cones}
The (polytopal) cone in $(\ker \pi)^*$ containing a point $\psi \in (\ker \pi)^*$ is the intersection
of all $\pi^\vee$-images of the cones in $\Delta (P)$ that contains $\psi$, that is,
$$\sigma(\psi):=[\psi] := \cap \{\pi^\vee N(P, F) | F \in L(P), \psi \in  \pi^\vee N(P, F)\}.$$
We  use $\Gamma^* = \Gamma^* (P, Q)$ to denote the set of all those cones, partially ordered by inclusion.
(These cones form the normal fan of the fiber polytope $\Sigma (P,Q)$ (\cite{BilleraSturmfels92})). 
\end{defn}

Consider the polytope $P_q = \pi^{-1} (q)$ for any $q \in Q$.
For any linear function $\psi \in (\ker \pi)^*$, let $P_q^\psi$ be the nonempty face of $P_q$
defined by $\psi$. The normal cone $N(P_q, P_q^\psi)=N(P_q^\psi)$
  depends only on the cell ${\bf c} = {\bf c} (q)$. 
So, we may use $\Delta({\bf c})$ to be the normal fan $\Delta (P_q) $ for all $q \in \relint {\bf c}$,
where ``relint'' denotes ``relative interior''. 
 For each $\psi$, let $[\psi]_{\bf c}$ denote the minimal cone in
$\Delta({\bf c})$ that contains $\psi$. 
Note that $\Delta({\bf c})$ is a common refinement of $\Delta({\bf c}')$ for all ${\bf c}' \le {\bf c}$.

In addition, for each face  
$P_q^\psi$ of $P_q$ there exists a unique minimal face of $P$ that contains $P_q^\psi$:
the intersection of all faces that contain $P_q^\psi$. We may use $F_{q, \psi}$ to denote this face
of $P$. Since it only depends on the cell ${\bf c} = {\bf c}(\psi)$,
  we may write it as $F_{{\bf c}, \psi}$.
Moreover, $F_{{\bf c}, \psi}$ only depends on the cone $[\psi]_{\bf c}$.
So, we oftentimes denote it by $F_{{\bf c},  [\psi]_{\bf c}}$.
Note that $F_{{\bf c}, \psi} = F_{{\bf c}, \psi'}$ if and only if $[\psi]_{\bf c}=[\psi']_{\bf c}$.
Note also that $N(P, F_{q, \psi})$ projects onto $N(P_q, P_q^\psi)$ under the projection $\pi^\vee$.

Having  the above explained, we can conclude that the fan of cones in the definition \ref{cones}
is the common refinement of
all the normal fans  $\Delta({\bf c})$ (${\bf c} \in \Gamma$).

\section{Coherent strings}
 
Given any linear function $\psi \in (\ker \pi)^*$, we have a collection of faces in $L(P)$
defined as follows:
$$\cF (\psi) = \{F_{{\bf c}, \psi} | {\bf c} \in \Gamma\}.$$

\begin{defn} (See \cite{Ziegler}) The above collection
$\cF (\psi) = \{F_{{\bf c}, \psi} |{\bf c} \in \Gamma \}$ of faces of $P$ is called a (global)
coherent string
of the projection $\pi: P \rightarrow Q$. We set $\cT_{coh}=\cT_{coh}(P,Q)$
 to be the set of all coherent strings of $\pi: P \rightarrow Q$, partially ordered by: 
$\cF (\psi) \le \cF (\psi') \; \hbox{if} \; \cup \cF (\psi) \subset \cup \cF (\psi').$
\end{defn}

\begin{prop} 
\label{coherentstring} Given any linear function $\psi \in (\ker \pi)^*$, the coherent
string $\cF (\psi)$ satisfies the following properties
\begin{enumerate}
\item $\{ \pi (F_{{\bf c}, \psi}) | {\bf c} \in \Gamma  \}$ 
       is  a polytopal subdivision of $Q$ without repetitions;
\item $\pi (F_{{\bf c}', \psi}) \subset \pi (F_{{\bf c}, \psi})$ if and only if 
         $F_{{\bf c}', \psi} =  \pi^{-1} (\pi(F_{{\bf c}', \psi})) \cap  F_{{\bf c}, \psi}$. 
\end{enumerate}
\end{prop}

\proof
Obviously, we have
$$\cup \{\pi (F_{{\bf c}, \psi}) | {\bf c} \in \Gamma\} = \cup {\bf c} = Q.$$
But, $q \in \relint \pi(F_{{\bf c}, \psi}) \cap \relint\pi(F_{{\bf c}', \psi}) \ne \emptyset$ if and only if 
$F_{{\bf c}, \psi} = \cap _{F \supset P^\psi_q} F = F_{{\bf c}', \psi}$. This proves (1).
(2) follows easily.
\endproof

\begin{lem}
\label{cohstr=coch}
$\cF (\psi)=\cF (\psi')$ if and only if $\sigma(\psi) = \sigma(\psi')$.
\end{lem}

\proof
$\cF (\psi)=\cF (\psi')$  if and only if $F_{{\bf c}, \psi}=F_{{\bf c}, \psi'}$ 
for all ${\bf c} \in \Gamma$
if and only if $[\psi]_{\bf c} = [\psi']_{\bf c}$ for all ${\bf c} \in \Gamma$ 
if and only if $\sigma(\psi) = \sigma(\psi')$.
\endproof

So, we may denote a coherent string by $\cF (\sigma)$ for some $\sigma \in \Gamma^*$.
Given any point $q \in Q$, we use $[q]_\sigma$ to denote the minimal face in $\cF (\sigma)$
such that its $\pi$-image contains $q$. Observe that
$F_{q, \sigma} = F_{q', \sigma}$ if and only if  $[q]_\sigma = [q']_\sigma$.

\begin{rem}
\label{moretoone}
Notice that two different coherent strings may induce identical subdivision
of $Q$. This will not be the case  if the projection is non-degenerate,
i.e., sends distinct vertices of $P$ to distinct points. In general, the issue
can typically be dealt with by keeping the label of each vertex.
\end{rem}

\begin{prop} (Billera-Sturmfels \cite{BilleraSturmfels92})
\label{cohstranticoch}
The poset $\cT_{coh}$ of all coherent strings and the poset $\Gamma^*$
are anti-isomorphic. That is, 
$\cT_{coh} \cong (\Gamma^*)^{opp}$ and $\Gamma^* \cong (\cT_{coh})^{opp}$.
\end{prop}

\proof Clearly the function $\cF: \Gamma^* \rightarrow \cT_{coh}$ is a bijection
by the definition and Lemma \ref{cohstr=coch}.
Let $\sigma$ be a cone and $\cF (\sigma)$ be its corresponding
coherent string. Assume that $\sigma \le \sigma$. Pick $\psi \in
\sigma$ and $\psi' \in \sigma$. Then $P_q^\psi \ge P_q^{\psi'}$.
 That is, $F_{q, \psi} \ge F_{q, \psi'}$. Since
$\cF (\sigma) = \{F_{q, \psi} | q \in Q\}$ and
$\cF (\sigma) = \{F_{q, \psi'} | q \in Q\}$,
we obtain that $\cF(\sigma) \ge \cF(\sigma)$. That is,
the function $\cF: \Gamma^* \rightarrow \cT_{coh}$ is an order-reversing bijection.
\endproof

\section{Locally coherent strings}

Locally coherent strings are natural generalization of coherent strings.
Geometrically, they are ``continuous liftings of polytopal subdivisions of $Q$''
to the face lattice of $P$.
It turns out that all of these things  can be constructed as follows.

Consider an arbitrary  map $$\widetilde{\Psi}: \Gamma \rightarrow (\ker \pi)^*.$$ 
It gives rise to a collection of faces of $P$ as follows
$$\cF(\widetilde{\Psi}) =\{ F_{{\bf c}, \widetilde{\Psi}({\bf c})} |{\bf c} \in \Gamma \}.$$

\begin{defn} (See \cite{Ziegler}) The collection $\cF(\widetilde{\Psi})$ of faces of $P$
is called   a locally coherent string if $\widetilde{\Psi}$ satisfies
the so-called {\it locally coherent condition}: that is, for any ${\bf c} \ge {\bf c}'$
$$ \relint [\widetilde{\Psi} ({\bf c})]_{\bf c} \subset 
\relint [\widetilde{\Psi} ({\bf c}')]_{{\bf c}'}.$$
We set $\cT =\cT (P, Q)$ to denote the set of all locally coherent strings,
partially ordered by:
$\widetilde{\Psi} \le \widetilde{\Psi'}$ if and only if 
$\cup \cF(\widetilde{\Psi}) \subset \cup \cF(\widetilde{\Psi'})$.

\end{defn}

One can check that the locally coherent condition is equivalent to
$$ \relint \pi^\vee (N(F_{{\bf c}, \widetilde{\Psi} ({\bf c})})) \subset
\relint \pi^\vee (N(F_{{\bf c}', \widetilde{\Psi} ({\bf c}')})).$$

\begin{prop} (\cite{BKS94}, \cite{RambauZiegler})
\label{stringtofaces}
Assume that a map $\widetilde{\Psi}: \Gamma \rightarrow (\ker \pi)^*$ satisfies
the locally coherent condition. Then the locally coherent string
 $\cF(\widetilde{\Psi}) $ enjoys
\begin{enumerate}
\item  $\{\pi(F)| F \in \cF(\widetilde{\Psi})\}$ is a subdivision of $Q$ without repetition.
\item  $\pi(F) \subset \pi(F')$ if and only if $F= F' \cap \pi^{-1}(\pi(F))$
    for any $F, F' \in  \cF(\widetilde{\Psi})$.
\end{enumerate}
Furthermore, every collection of the faces of $P$ satisfying (1) and (2) arises in this way.
\end{prop}

\proof The theorem was  proved in \cite{BKS94} for the case when $\dim Q=1$ and
was explicitly pointed out in \cite{RambauZiegler} for the general cases.

A key observation in the proof is that when ${\bf c}' <{\bf c}$, the locally coherent condition
implies that $F_{{\bf c}', \widetilde{\Psi}({\bf c}')}$ is either identical to
$F_{{\bf c}, \widetilde{\Psi}({\bf c})}$ or is a proper face of $F_{{\bf c}, \widetilde{\Psi}({\bf c})}$.
(1) and (2) basically follow from this observation.

Conversely, given any collection $\cF$ satisfying (1) and (2). For any ${\bf c} \in \Gamma$,
let $F({\bf c})$ be the unique minimal face in $\cF$ such that
$\relint {\bf c} \subset \relint \pi(F({\bf c}))$. Let 
$\psi ({\bf c}) \in \relint \pi^\vee(N(F({\bf c}))) \subset
(\ker \pi)^*$ such that $F({\bf c}) = F_{{\bf c}, \psi ({\bf c})}$. Define 
$$\widetilde{\Psi}: \Gamma \rightarrow (\ker \pi)^*$$
as
$$\widetilde{\Psi}: {\bf c} \rightarrow \psi ({\bf c}).$$
Now, if ${\bf c}' \le {\bf c}$ (we may well assume that
${\bf c}$ covers ${\bf c}'$), either $F({\bf c}') = F ({\bf c})$ (in this case,
$\relint {\bf c}' \subset \relint \pi(F({\bf c}))$) or $F({\bf c}')$ is a proper face of
$F({\bf c})$ (in this case, ${\bf c}'$ lies on the boundary of $\pi(F ({\bf c}))$).
In the first case, we have $[\psi ({\bf c})]_{\bf c}$ equals $[\psi ({\bf c}')]_{{\bf c}'}$ which obviously implies that  
 $\relint [{\psi} ({\bf c})]_{{\bf c}} \subset \relint [{\psi} ({\bf c}')]_{ {\bf c}'}$.
Assume we are in the second case. We then have that
$[\psi ({\bf c})]_{\bf c}$ divides $[\psi ({\bf c}')]_{{\bf c}'}$ from its relative interior which also implies that  
 $\relint [{\psi} ({\bf c})]_{{\bf c}} \subset \relint [{\psi} ({\bf c}')]_{ {\bf c}'}$.

This completes the proof.
\endproof




The minimal elements of $\cT$ tie to the following

\begin{defn} $\cF(\widetilde{\Psi})$ is said to be tight if 
$$\dim \pi(F) = \dim F$$ for all $F \in
\cF(\widetilde{\Psi})$. (That is, $F$ does not drop dimension under the projection.)
\end{defn}

\begin{prop}
\label{tight=minimal}
$\cF (\widetilde{\Psi})$ is minimal if and only if $\cF(\widetilde{\Psi})$ is tight.
\end{prop}

\proof
This is  Lemma 9.5 of \cite{Ziegler} whose proof for the coherent case extends to
the general case. The proof goes as follows. Assume that $F(\widetilde \Psi)$ is minimal
but there exists a face $F \in F(\widetilde \Psi)$ such that the projection $F \to \pi(F)$
drops dimension. Now consider the polytopal projection $F \to \pi(F)$,
 take any non-trivial coherent string $\Psi_F$ of this projection, and substitute
$F$ from $\cF(\widetilde \Psi)$ by $\cF(\Psi_F)$, we will get a new locally coherent
string of the original projection which is certainly smaller than $\cF(\widetilde \Psi)$.
Hence $\cF(\widetilde \Psi)$ has to be tight.

The converse is easy. We omit further details.
\endproof

\section{Coherent costrings}

Let $q \in Q$ be any interior point. The point $q$ selects the following
collection of cones in $\Delta (P)$
$$\Delta (q) = \{ N(P, F_{q, \psi}) | \psi \in (\ker \pi)^* \}.$$
The collection only depends on the cell ${\bf c} (q)$. That is,

\begin{lem}
\label{bj}
$\Delta (q) = \Delta (q')$ if and only if ${\bf c}(q) ={\bf c}(q')$.
\end{lem}

\proof
$\Delta (q) = \Delta (q')$  if and only if $F_{q, \sigma}=F_{q', \sigma}$ 
for all $\sigma \in \Gamma^*$
if and only if $[q]_{\sigma} = [q']_{\sigma}$ for all $\sigma \in \Gamma^*$ 
if and only if ${\bf c}(q) ={\bf c}(q')$.
\endproof

Thus
we may denote $\Delta (q)$ by $\Delta ({\bf c})$ if $q \in \relint {\bf c}$ (${\bf c} \in \Gamma$). 

\begin{defn}  $\Delta ({\bf c})$ is called a coherent
costring. Let $\cT^*_{coh} = \{ \Delta ({\bf c}) | {\bf c} \in \Gamma\}$ 
be the set of all coherent costrings, partially ordered by:
$$\Delta ({\bf c}) \le \Delta ({\bf c}') \; \hbox{if}\; \cup \Delta ({\bf c}) \subset  \cup \Delta ({\bf c}').$$ 
\end{defn}

\begin{prop} Given any ${\bf c} \in \Gamma$. The collection $\Delta ({\bf c})$ satisfies the following
properties
\begin{enumerate}
\item $\{\pi^\vee (N(P, F_{{\bf c}, \psi})) | \psi \in (\ker \pi)^*  \}$ 
       is  a cone subdivision of $(\ker \pi)^*$ without repetitions.
\item $\pi^\vee (N(P, F_{{\bf c}, \psi})) \subset \pi^\vee (N(P, F_{{\bf c}', \psi}))$
      if and only if $$N(P, F_{{\bf c}, \psi})= (\pi^\vee)^{-1} (\pi^\vee (N(P, F_{{\bf c}, \psi}))) \cap
   \pi^\vee (N(P, F_{{\bf c}', \psi})).$$    
\end{enumerate}
Moreover, $\{\pi^\vee(\Delta(q))\}$ equals the normal fan $\Delta(P_q)$ if $q \notin \partial Q$.
\end{prop}

\proof The proof of (1) and (2) is analogous to that of Proposition \ref{coherentstring} and is a careful check using definitions. In section \ref{lcc}, we will prove a general version (i.e., for locally coherent costrings)
which includes
this as a particular case.
For the last statement, just observe that $\pi^\vee (N(P, F_{q,\psi})) = N (P_q, P_q^\psi)$.
We omit further details.
\endproof

\begin{rem}
Like in Remark \ref{moretoone}, two different coherent costrings may induce identical
cone subdivision of $(\ker \pi)^*$. 
\end{rem}

\begin{prop}
The poset $\cT^*_{coh}$ of all coherent costrings and the poset $\Gamma$
of all cells are anti-isomorphic. That is, 
$\cT^*_{coh} = \Gamma^{opp}$ and $\Gamma = (\cT^*_{coh})^{opp}$.
\end{prop}

\proof Treat $\Delta$ as a function from $\Gamma$ to $\cT^*_{coh}$.
This is then a bijection by Lemma \ref{bj}.
Assume that ${\bf c}' \le {\bf c}$ are two cells.
Pick $q_{\bf c}' \in \relint {\bf c}'$ and $q_{\bf c} \in \relint {\bf c}$.
We have that $F_{q_{\bf c}', \psi} \le F_{q_{\bf c}, \psi}$ for all
$\psi \in (\ker \pi)^*$. That is,
$N(P,F_{q_{\bf c}', \psi}) \ge N(P, F_{q_{\bf c}, \psi})$ for all $\psi$.
By definition, we obtain $\Delta({\bf c}') \ge \Delta({\bf c})$.
That is, $\Delta$ is an  order-reversing bijection.
\endproof

\section{Locally coherent costrings}
\label{lcc}

Locally coherent costrings are generalizations of (global) coherent costrings. Geometrically,
they are $\pi^\vee$-liftings of conical subdivisions of $(\ker p)^*$.

\begin{defn}
\label{lccostring} Let $\Delta_0$ be an {\it arbitrary} fan in $V^*$.
A subcollection $\Delta$ of cones in $\Delta_0$ is called a locally coherent costring of
 $\pi^\vee: |\Delta_0| \subset V^* \rightarrow \pi^\vee (|\Delta_0|) \subset (\ker \pi)^*$ 
if it satisfies that
\begin{enumerate}
\item $\{\pi^\vee(\sigma) | \sigma \in \Delta \}$ form a cone subdivision of $ \pi^\vee(|\Delta_0|)$ without
       repetition.
\item  $\pi^\vee(\sigma) \subset \pi^\vee(\sigma')$ 
if and only if $\sigma = (\pi^\vee)^{-1}(\pi^\vee (\sigma)) \cap \sigma'$.
\end{enumerate}
\end{defn}

It follows from the definition that if $\sigma \ne \sigma' \in \Delta$, then
$\pi^\vee (\sigma)$ and $\pi^\vee (\sigma')$ must intersect along a proper face. In particular,
$\pi^\vee (\Delta) = \{\pi^\vee (\sigma) | \sigma \in \Delta\}$ is a {\it fan} in $W$. 
(Warning: a fan in our sense may not be strongly convex.)

The following picture illustrates an example of a (locally) coherent
costring $\Delta$ in ${\Bbb R}^3$ with the projection
$\pi^\vee (\Delta)$ in ${\Bbb R}^2$. (See \cite{RambauZiegler} for pictures of locally coherent strings.)

\vskip   1.9in
\begin{center}
Figure 1
\end{center}

The most interesting case where rich structures and duality
emerge is the case when $\Delta_0$ is the normal fan $\Delta (P)$ of the polytope $P$.

\begin{defn}
The set of locally coherent costrings of $\pi^\vee: |\Delta (P)| \rightarrow (\ker \pi)^*$
is denoted by $\cT^*= \cT^*(P,Q)$. This set is equipped
with a natural partial order by inclusion. Note that 
$\Gamma^*$ is a common refinement of $\pi^\vee (\Delta)$ for all $\Delta \in \cT^*$.
\end{defn}

As before, we will actually construct all locally coherent costrings.
So,  again we consider a map $$\widetilde{\Delta}: \Gamma^* \rightarrow Q$$
satisfying the locally coherent condition: for any $\sigma > \sigma'$,
$$\relint [\widetilde{\Delta}(\sigma)]_\sigma \subset \relint [\widetilde{\Delta}(\sigma')]_{\sigma'}.$$
  $\widetilde{\Delta}$ gives rise to a 
collection $\cC(\widetilde{\Delta})$ of cones in $\Delta(P)$
as follows
$$\cC(\widetilde{\Delta}) =\{N(P, F_{\widetilde{\Delta}(\sigma),\sigma}) | \sigma \in \Gamma^* \}.$$

\begin{thm}
Let $\widetilde{\Delta}: \Gamma^* \rightarrow Q$ satisfy the locally coherent
condition.  $\cC(\widetilde{\Delta})$
is a locally coherent costring of $\pi^\vee: |\Delta (P)| \rightarrow (\ker \pi)^*$.
Furthermore, every locally coherent costring of $\pi^\vee: |\Delta (P)| \rightarrow (\ker \pi)^*$
arises in this way.
\end{thm}

\proof The proof of this theorem is in a dual manner  analogous to the case of locally coherent strings. 

We see immediately that
$$\cup_{\sigma \in \Gamma^*} \pi^\vee (N(F_{\widetilde{\Delta}(\sigma), \sigma})) \supset
  \cup_{\sigma \in \Gamma^*} \sigma = (\ker \pi)^*$$
which implies the equality. 
To see that we have a subdivision without repetition, just observe that if $\sigma > \sigma'$, either 
$F_{\widetilde{\Delta}(\sigma), \sigma} = F_{\widetilde{\Delta}(\sigma'), \sigma'})$
or $F_{\widetilde{\Delta}(\sigma), \sigma}$ is a face of $F_{\widetilde{\Delta}(\sigma'), \sigma'})$
by the locally coherent condition.
Condition (2) of Definition \ref{lccostring} follows from the same observation.

Conversely, for any $\sigma \in \Gamma^*$, let $N(F(\sigma))$ be the minimal cone in 
$\cC $
such that $\pi^\vee (N(F(\sigma)))$ contains $\sigma$. Let $q(\sigma) \in \relint F(\sigma)$ be a point such that
$F(\sigma) = F_{q(\sigma), \sigma}$. Define
$$\widetilde{\Delta}: \Gamma^* \rightarrow Q$$ by
$$\widetilde{\Delta}: \sigma \rightarrow q(\sigma).$$
If $\sigma > \sigma'$
(it suffices to assume that $\sigma$ covers $\sigma'$, that is,
$\dim \sigma = \dim \sigma' + 1$), since the subdivision of $Q$ induced by the coherent string defined by $\sigma$
refines the subdivision of $Q$ induced by the coherent string defined by $\sigma'$,  we have either
$F(\sigma) = F(\sigma')$ or $F(\sigma)$ divides $F(\sigma')$ from its relative interiors. In either case,
it implies that $\relint [q(\sigma)]_\sigma \subset \relint [q(\sigma')]_{\sigma'}$.
In other words,
the map $\widetilde{\Delta}$ indeed satisfies the locally coherent condition.
\endproof



\begin{defn}
A locally coherent costring $\cC(\widetilde{\Delta})$ is tight if 
$\dim \pi^\vee (\sigma) = \dim \sigma$ for every $\sigma \in \cC(\widetilde{\Delta})$.
\end{defn}

Analogous to Proposition \ref{tight=minimal} for locally coherent strings, 
we have the following proposition
 for locally coherent costrings whose proof is similar 
to that of Proposition \ref{tight=minimal} and is left to the reader.

\begin{prop}
A locally coherent costring $\cC(\widetilde{\Delta})$ is tight if and only if it is minimal.
\end{prop}

\begin{rem}
In view that $\Gamma^* = \cT_{coh}$, every locally coherent costring $\cC(\Delta)$ 
singles out a unique minimal element from any given  coherent string $ \cF (\sigma) \in \cT_{coh}$:
 i.e., the minimal face in $\cF (\sigma)$ such that its image contains the cell 
$\widetilde{\Delta}(\sigma)$.
\end{rem}

\section{Virtual chambers}

{\sl Virtual cells} are generalizations of the real cells.
In toric geometry, real cells define 
{\it projective}  quotient varieties, while
virtual cells define (not-necessarily projective) quotient varieties.

We refer to \cite{dHSS} for a treatment  of virtual chambers in the case of  triangulations of a
point set configuration. 

Any map $\widetilde{{\bf c}}: \Gamma^* \rightarrow Q$ gives rise to a
collection of faces of $P$ as follows
$$\cF(\widetilde{{\bf c}}) = \{F_{{\widetilde{{\bf c}} (\sigma), \sigma}} | \sigma \in \Gamma^*\}.$$

\begin{defn} $\cF(\widetilde{{\bf c}})$ is called a virtual cell if $\widetilde{{\bf c}}$ satisfies
the following  locally coherent condition: for $\sigma \ge \sigma'$
$$ \relint [\widetilde{{\bf c}} (\sigma)]_\sigma 
  \subset \relint [\widetilde{{\bf c}} (\sigma')]_{\sigma'}.$$
We use $\Gamma_{vir} = \Gamma_{vir}(P, Q)$ to denote the set of all  virtual cells,
 partially ordered by inclusion:
$\widetilde{{\bf c}} \le \widetilde{{\bf c}}'$ if $\cup \cF (\widetilde{{\bf c}}) \subset 
\cup \cF (\widetilde{{\bf c}}')$. Maximal virtual cells are called virtual chambers.
\end{defn}

\begin{thm} 
\label{chamberstofaces}
Assume that a map
$\widetilde{{\bf c}}: \Gamma^* \rightarrow Q$ satisfies the locally coherent condition.
Then the virtual cell $\cF(\widetilde{{\bf c}})$  satisfies
the following property: for any coherent string $\Psi$, $\cC (\widetilde{{\bf c}}) \cap \cF(\Psi)$ is nonempty and         contains exactly one minimal element.
Furthermore, every collection of the faces of $P$ satisfying  the above property arises in this way.
\end{thm}

\proof
Let $\cF = \cF (\sigma) \in \cT_{coh}$ with $\sigma \in \Gamma^*$.
Clearly both $\cF (\sigma) $ and $\cF(\widetilde{{\bf c}})$ contain
the element $F_{\widetilde{{\bf c}} (\sigma), \sigma}$. That is, the intersection is not empty.
Now given any two faces $F_{p, \sigma}$ and $F_{q, \sigma}\in \cF (\sigma) \cap \cF(\widetilde{{\bf c}})$,
either they have a common face or they do not meet at all. If they do not meet at all,
$\pi^\vee (N(F_{p, \sigma}))$ and $\pi^\vee (N(F_{q, \sigma}))$ 
will not have a non-trivial common face. Hence they meet
in relative interiors (both contain $\sigma$). But as members of the locally coherent costring
$\cC (\Delta)$ (transported
from $\Delta = \widetilde{c}$), they can not meet in relative interiors.
This shows that $F_{p, \sigma}\cap F_{q, \sigma}$ is a common face of each. Clearly, 
$\cF (\sigma) \cap \cF(\widetilde{{\bf c}})$ must,
as a subset of
a coherent string $\cF (\sigma)$ and with property that every two members have a common face,
have a unique minimal element.

On the other hand,  given a collection $\cF$ that
meets every coherent string $\cF (\sigma)$ in  exactly one minimal element.
Let $F(\sigma) \in \cF \cap \cF (\sigma)$ be this face. Pick a point
$q_\sigma \in \relint \pi(F(\sigma))$. We obtain a map
$$\widetilde{\bf c}: \Gamma^* \rightarrow Q$$ 
by sending $\sigma$ to $q_\sigma$.
We now check that
this map satisfies the locally coherent condition.
For any $\sigma' \le \sigma$, we have that the subdivision 
$\{\pi(F) | F \in \cF (\sigma) \}$ refines the subdivision 
$\{\pi(F) | F \in \cF (\sigma') \}$. It follows that
$\relint \pi(F(\sigma))  \subset  \relint \pi(F(\sigma'))$
because either $F(\sigma) = F(\sigma')$ or $\pi(F(\sigma))$ divides
$\pi(F(\sigma))$ in its relative interiors.
 This checks the locally coherent condition.
\endproof

\begin{prop} The following two are equivalent. 
\begin{enumerate}
\item  For any coherent string $\Psi$, $\cF (\widetilde{{\bf c}}) \cap \cF(\Psi)$ is nonempty and
            contains exactly one minimal element.
\item For any locally coherent string $\widetilde{\Psi}$, 
            $\cF (\widetilde{{\bf c}}) \cap \cF(\widetilde{\Psi})$ is nonempty and
            contains exactly one minimal element.
\end{enumerate}
\end{prop}
\proof
For point configurations in general position, this is  proved in \cite{dHSS} for virtual chambers.

(2) obviously implies (1). To see the other direction, assume $\cF (\widetilde{{\bf c}})$
satisfies (1) but not (2), then there will be two faces $F$ and $F'$ 
in some locally coherent string $\widetilde{\Psi}$ such that they are both contained in 
$\cF (\widetilde{{\bf c}})$ and are two distinct minimals in
$\cF (\widetilde{{\bf c}}) \cap \cF(\widetilde{\Psi})$. 
Since $\pi (F) \cap \pi (F')$ must be  empty, 
we can find a coherent string $\Psi$ to contain both $F$ and $F'$. This contradicts to that
$\cF (\widetilde{{\bf c}})$
satisfies (1).
\endproof

This equivalence will not, however,  be used later.

\section{Virtual cones}

Virtual cones are simply dual versions of virtual cells.

An arbitrary  map $\widetilde\sigma: \Gamma \rightarrow (\ker \pi)^*$
defines a collection $\cC (\widetilde\sigma)$ of cones of $\Delta(P)$
$$\cC (\widetilde\sigma)=\{N(F_{{\bf c}, \widetilde\sigma ({\bf c})}, P)| {\bf c} \in \Gamma \}.$$

\begin{defn} $\cC (\widetilde\sigma)$ of cones of $\Delta(P)$
is called a virtual cone if $\widetilde\sigma$ satisfies the locally coherent condition. 
The set of all virtual cones is denoted by $\Gamma^*_{vir} = \Gamma^*_{vir} (P,Q)$,
 partially ordered by inclusion.
\end{defn}

\begin{thm}
Let $\widetilde\sigma: \Gamma \rightarrow (\ker \pi)^*$ satisfy the locally coherent
condition. Then the  virtual cone $\cC (\widetilde\sigma)$
satisfies the following property.
For every coherent costring $\cC({\bf c})$, $\cC (\widetilde\sigma) \cap \cC({\bf c})$
is non-empty and contains exactly one maximal member.
Furthermore, every collection of cones of $\Delta(P)$ that satisfies above property
arises in this way.
\end{thm}

\proof Transporting $\cC({\bf c})$ and $\cC (\widetilde\sigma)$
 to their corresponding face collections $\cF ({\bf c})$
 and $\cF (\widetilde \Psi)$ 
where $\widetilde \Psi = \widetilde\sigma$, we see immediately that
$\cF (\widetilde \Psi) \cap \cF ({\bf c})$ contains a unique minimal element (i.e., the intersection
of all the faces of $\cF (\widetilde \Psi)$ that contain ${\bf c}$).
Translating back, this means that  $\cC (\widetilde\sigma) \cap \cC({\bf c})$
contains exactly one maximal member.
\endproof



\begin{rem} Again, the following should  be equivalent.
\begin{enumerate}
\item $\cC (\widetilde{\sigma})$ meets every coherent costring 
       in exactly one maximal member.
\item $\cC (\widetilde{\sigma})$ meets every locally coherent costring 
       in exactly one maximal member.
\end{enumerate}
\end{rem}

\section{Summary of dualities}

It may be worthwhile to summarize some of the results in the previous section as the following dualities.

\begin{thm} The following always holds.
\begin{enumerate}
\item A collection ${\cF}$ of the faces of $P$ is a 
      locally coherent string of $\pi:P \rightarrow Q$
     if and only if $\cC = \{N(P,F)|F \in \cF\}$ is a 
       virtual cone of $\Delta (P) \rightarrow (\ker \pi)^*$.
\item A collection ${\cF}$ of the faces of $P$ is a  virtual
      cell of $\pi:P \rightarrow Q$
      if and only if $\cC = \{N(P,F)|F \in \cF\}$ is a 
      locally coherent costring of $\Delta (P) \rightarrow (\ker \pi)^*$.
\item A collection ${\cC}$ of the cones  of $\Delta(P)$ is a 
      locally coherent costring of $\Delta(P) \rightarrow (\ker \pi)^*$
     if and only if $\cF = \{F| N(P,F) \in \cC\}$ is a 
       virtual cell of $P \rightarrow Q$.
\item A collection ${\cC}$ of the cones  of $\Delta(P)$ is a virtual
      cone of $\Delta (P) \rightarrow (\ker \pi)^*$
      if and only if $\cF = \{F|N(P,F) \in \cC\}$ is a 
      locally coherent string of $\pi:P \rightarrow Q$.
\end{enumerate}
Furthermore, the correspondence in each of (1), (2), (3), (4)
reverses the inclusion-induced partial orders.
\end{thm}

\proof
This theorem is a corollary to the combination of the results in the previous sections.
\endproof

Let $\Mor'(S, S')$ be the poset of order-reversing maps between two posets $S$ and $S'$.
The we have
\begin{enumerate}
\item $\cT \subset \Mor'(\Gamma, \Gamma^*)$ and elements in $\cT_{coh}$ correspond to constant maps;
\item $\cT^* \subset \Mor'(\Gamma^*, \Gamma)$ and elements in $\cT_{coh}$ correspond to constant maps;
\item $\Gamma_{vir}^* \subset \Mor'(\Gamma, \Gamma^*)= \Mor' (\cT^*_{coh}, \cT_{coh})$ and elements in $\Gamma^*$ correspond to constant maps;
\item $\Gamma_{vir} \subset \Mor'(\Gamma^*, \Gamma) = \Mor' (\cT_{coh}, \cT_{coh}^*)$ and elements in $\Gamma$ correspond to constant maps;
\end{enumerate}
Note the natural correspondences amongst the bigger posets.
These relations seem to worth some investigation. In particular, it would be interesting to know
if an arbitrary element in each of the bigger posets corresponds to some meaningful combinatorial
construction.



\section{Realizable virtual chamber.}
\label{deformcells}

Examples indicate that
 most of virtual cells can be obtained by topologically deforming  the original polytopal cell complex.
We now explain this.

Let $P_t$ be a topological deformation of $P$,
while the lattice structure on $L(P)$
or all of the face relations in $L(P)$ are memorized. Then the image $Q_t =\pi (P_t)$
will be called a tamed deformation of $Q = Q_0$. Note that the deformation also induces
a (tamed) deformation of walls, and in particular, a (tamed) deformation $\Gamma_t$ of the chamber
complex $\Gamma$.

Another possible way to think of a tamed deformation is to interpret  it as a deformation $\pi_t$
of the affine projection $\pi = \pi_0$. In that way, all the twisted projections $\pi_t$ ($t\ne 0$)
may be called  {\it virtual affine projections}. Thus a tamed deformation $Q_t$ and its induced chamber complex
$\Gamma_t$ will just be the image of $P$ under the virtual affine projection $\pi_t$.

\begin{exmp} 
\label{deformcell}
Figure 2 is a typical tamed deformation of a hexagon and its internal wall structures.
Note that by the deformation, we gain one more (deformed) chamber, three more deformed 1-dimensional cells,
and two more 0-dimensional cells (gain three but loose one).

The next picture is not a tamed deformation because an interior wall and a face of the hexagon
intersect at a wrong place. (Notice the gain of an extra chamber in the upper right corner.)
\end{exmp}

\vskip   1.8in
\begin{center}
Figure 2
\end{center}

\vskip   1.8in
\begin{center}
Figure 3
\end{center}

\begin{rem}
In terms of geometry of the torus action, the above means that one can not alert
the intersection relations
among fixed point varieties (for various subtori)
under the deformation.
\end{rem}

\begin{prop}
Every tamed deformation of $P$ induces a tamed deformation of the $\pi$-induced
polytopal subdivision $\pi(T)$ of $Q$ for any given $T \in \cT$. In particular, the combinatorial type of
$\pi(T)$ remains the same under the deformation.
\end{prop}

It follows that if $\Gamma_t$ is the cell complex induced by $P_t$, then
$\Gamma_t$ is the common refinement of $\pi(T_t)$ for all $T \in \cT$.

\begin{cor}
A deformed cell in $\Gamma_t$ corresponds canonically to
a virtual cell in $\Gamma$.
\end{cor}
\proof
Each cell of $\Gamma_t$ specifies a unique minimal deformed subpolytope in $\pi(T_t)$ (thus
a unique  minimal face $F_t$ in $T_t$) for every $T \in T$.
 This leads to a unique minimal face $F$ in the original $T$. 
By Theorem \ref{chamberstofaces},
we conclude that  every cell of $\Gamma_t$ gives rise
to a virtual cell of $\pi: P \rightarrow Q$.
\endproof

\begin{defn}
A virtual cell (chamber) that corresponds to some deformed cell (chamber) is called a
{\sl realizable} virtual cell (chamber).
\end{defn}

\begin{defn}
A deformed wall in $Q_t$ will be called a virtual wall. Of Figure 2, the three bended walls
are typical virtual walls.
\end{defn}

Some realizable virtual chambers deserve special attention.

\begin{defn} (\cite{dHSS})
A cell in $\Gamma_t$ ($t=0$ or otherwise) that contains an extremal vertex
is  called  lexicographic. 
\end{defn}

It is easy to see that the realizable virtual cell that corresponds to a lexicographic cell
is independent of the deformation. From this, we get

\begin{cor}
The subpost of realizable virtual cells is connected.
\end{cor}

\begin{rem}
It is not clear if all virtual cells are realizable. By the work of \cite{BKS94},
it looks that it is the case when $\dim Q =1$ and perhaps also when $\dim Q =2$.
In general, results on the generalized Baues conjecture suggests that it should fail in general.
\end{rem}


\section{Toric varieties}

Let $X^d$ be an arbitrary toric variety over the field of complex numbers, not necessarily quasi-projective,
acted upon by the open dense torus $({\Bbb C}^*)^d$. 
We will study categorical (good) quotients of $X$ by a subgroup $G$ of $({\Bbb C}^*)^d$. 
The method renders it useful the combinatorial concepts of \S 2.

Fix a lattice $M \cong {\Bbb Z}^d$ and its dual lattice $N$. A fan
$\Delta_0$ in $N_{\Bbb R} = N \otimes {\Bbb R}$ consists of a finite collection
of strongly convex rational polyhedral cones 
$\sigma \subset N_{\Bbb R}$ which is closed under intersection and 
taking faces.

Given (and fix) the  fan  $\Delta_0$  and a cone $\sigma$ in the fan,
the dual cone 
$$\sigma^\vee =\{ v \in M_{\Bbb R}: \langle v, \sigma \rangle \ge 0 \}$$
determines the semigroup algebra ${\Bbb C}[\sigma^\vee \cap M]$ which is finitely generated
and the affine variety $$A_{\sigma} = \Spec ({\Bbb C}[\sigma^\vee \cap M]).$$
It turns out that these affine varieties can be glued together to get a toric variety
$$X_{\Delta_0}= \bigcup_{\sigma \in \Delta_0} \Spec ({\Bbb C}[\sigma^\vee \cap M]).$$

Recall that  $T_N = N \otimes {\Bbb C}^* = \Hom_{\Bbb Z} (M, {\Bbb C}^*)$ 
acts on ${\Bbb C}[M]$ by the formula
$$ \gamma_n (t) \chi^m = t^{\langle n, m \rangle} \chi^m$$
where $\gamma_n$ is the one-parametric subgroup generated by $n \in N$.

$X_{\Delta_0}$ can also be constructed as a quotient of certain open subset of some affine space
by a diagonalizable group. We refer the reader to \cite{Cox95} or \cite{Cox96} for more details.

\section{Locally coherent costrings and quotients}

Let $M^2$ be any sublattice of $M$. Consider the exact sequence
$$0 \rightarrow  M^2 \rightarrow  M @>{\pi}>>  M^1 \rightarrow  0.$$
If we apply $\Hom_{\Bbb Z} (-, {\Bbb C}^*)$ to it, then we get the exact sequence
$$1  \rightarrow \Hom_{\Bbb Z} (M^1, {\Bbb C}^*) \rightarrow \Hom_{\Bbb Z} (M, {\Bbb C}^*) \rightarrow \Hom_{\Bbb Z} (M^2, {\Bbb C}^*) \rightarrow 1.$$
The group $G = \Hom_{\Bbb Z} (M^1, {\Bbb C}^*)$ is the product of the torus 
$\Hom_{\Bbb Z} (M^1_{free}, {\Bbb C}^*)$ and the finite group
$\Hom_{\Bbb Z} (M^1_{tor}, {\Bbb C}^*)$.
Let $$0 \rightarrow  N^1 \rightarrow  N @>{\pi^\vee}>> N^2 \rightarrow  0$$
be the dual sequence.
(We point out that we also use $\pi$ ($\pi^\vee$) to denote its ${\Bbb R}$-linear extension.)

We need a couple of lemmas to proceed

\begin{lem} 
\label{lem:basic1}
Let $\sigma_0 \in \Delta_0$ be such that $\pi^\vee (\sigma_0) = \sigma$ is strongly convex. Then
\begin{enumerate}
\item ${\Bbb C}[\sigma_0^\vee \cap M]^G \cong {\Bbb C}[\sigma^\vee \cap M^2]$;
\item $A_{\sigma_0}/\!/G$ exists and is naturally identified with $A_{\sigma}$.
\end{enumerate}
\end{lem}

\proof 
Assertion (2) obviously follows from (1).
We only need to prove (1). On one hand, it follows from the definitions that 
${\Bbb C}[\sigma^\vee \cap M^2] \subset {\Bbb C}[\sigma_0^\vee \cap M]^G$. On the other hand,
any point of ${\Bbb C}[M]$ that is left invariant under the action of $G=\Hom_{\Bbb Z} 
(M^1, {\Bbb C}^*)$
must belong to $M_2$ (one can check this by the action formula given in \S 11). 
This leads to $[\sigma_0^\vee \cap M]^G \subset \sigma^\vee \cap M^2 $.
The lemma thus follows.
\endproof

The lemma below gives a criterion for the separateness of a quotient

\begin{lem}
\label{lem:basic2}
Let $\sigma$ and $\sigma'$ be two cones in $\Delta_0$.
The open subset $A_{\sigma} \cup A_{\sigma'}$ has the separated
categorical quotient $(A_{\sigma} \cup A_{\sigma'})/\!/G$ if and only if 
$\pi^\vee (\sigma) \cap  \pi^\vee (\sigma') $ is a face of each.
\end{lem}

\proof
By  Lemma \ref{lem:basic1}, the (possibly non-separated)
quotient variety 
$ A_{\sigma} \cup A_{\sigma'}/\!/G$ is naturally identified with
the variety $A_{\pi^\vee (\sigma)} \cup A_{\pi^\vee (\sigma')}$.
By a standard fact from the theory of toric varieties (see Lemma 1.4 of \cite{Fulton93}),
this is separated if and only if $ \pi^\vee (\sigma) \cap  \pi^\vee (\sigma')$ is a common face of both.
\endproof

\begin{defn}
A $G$-equivariant map $X \rightarrow Y$ is called a good quotient if $X$ admits a covering
by invariant affine open subsets $U_i = \Spec (S_i)$ such that $Y = \cup_i U_i/\!/G$
where $U_i/\!/G = \Spec (S_i^G)$.
\end{defn}


\begin{defn} A locally coherent costring of $\pi^\vee: N \rightarrow  N^2$
 is called {\it strongly convex} if every cone $\pi^\vee (\sigma)$ ($\sigma \in \Delta$) is strongly convex. 
\end{defn}

Recall that a locally coherent costring is tight (minimal) if the dimension of any cone
in the locally coherent costring does not drop under the projection.

Here comes the main theorem of this section.

\begin{thm}
\label{mainquotientthm}
Let $\Delta$ be a strongly convex locally coherent costing of $N \rightarrow  N^2$.
Then $U(\Delta)= \cup_{\sigma \in \Delta} A_\sigma \subset X$ is an invariant Zariski open subset
such that
\begin{enumerate}
\item     the separated categorical good quotient $U(\Delta) \rightarrow U(\Delta)/\!/G$ exists.
\item     $ U(\Delta)/\!/G$ is the toric variety
         $X_{\pi^\vee(\Delta)}$ defined by the induced fan $\pi^\vee (\Delta)$.
          \item    $U(\Delta) \rightarrow U(\Delta)/\!/G$ is a geometric quotient if and
       only if $\Delta$ is tight (minimal). 
\end{enumerate}
\end{thm}

\proof That $U(\Delta)$ is an invariant Zariski open subset is obvious.

(1) and (2) follow from the combination of Lemmas \ref{lem:basic1} and \ref{lem:basic2} and the
construction of toric varieties. For much of it, it suffices to observe that
$U(\Delta)/\!/G =\cup_{\sigma \in \Delta} A_\sigma/\!/G =
\cup_{\sigma \in \Delta} A_{\pi^\vee (\sigma)}$ (by Lemma 12.1 (2)).
Statement (3) is local. So to show it, it suffices to prove that $A_\sigma/\!/G$ is a geometric quotient for each
$\sigma \in \Delta$. $A_\sigma/\!/G$ being geometric is equivalent to that for every $x \in A_\sigma$,
$G_x$ is finite. Let $N_\sigma$ be the lattice generated by $\sigma$. Then by \S 3.1 (or \S 2.1) of Fulton \cite{Fulton93},
$T_{N_\sigma} = N_\sigma \otimes {\Bbb C}^*$ is the identity component
of the  isotropy subgroup of any point in the $T$-orbit ${\cal O}_\sigma \subset A_\sigma$ determined by
$\sigma$. Moreover, the identity component
of the  isotropy subgroup of any point in $A_\sigma$ is contained in $T_{N_\sigma}$.
So a point in $A_\sigma$ having isotropy subgroup in $G = T_{N^1}$ of positive dimension
if and only if $\dim T_{N^1} \cap T_{N_\sigma} >0$ if and only if
$N^1 \cap N_\sigma \ne \{0\}$ if and only if $\sigma$ drops dimension under the projection to $N^2$.
 This proves the theorem.
\endproof


As an immediate corollary of Theorem \ref{mainquotientthm} (3), we have
\begin{cor} 
If $\Delta$ is tight, then $U(\Delta)$ is the toric variety which is a principal
$G$-bundle over $U(\Delta)/\!/G$, assuming that $G$ acts on $U(\Delta)$ freely.
Moreover, all principal $G$-bundle among toric varieties arise this way.
\end{cor}

\begin{rem} As in the very general case of a quotient theory, degenerate quotients exist.
These are the quotients whose dimensions are less than expected (that is, smaller than $\dim X - \dim G$).
In our toric situation, this will be the case if the induced fan $\pi^\vee(\Delta)$ is not
strongly convex (that is, some projection $\pi^\vee(\sigma)$ contains a non-trivial vector space,
or equivalently, $\sigma \cap N^1$ is infinite). In this case,
$ U(\Delta) /\!/G$ is the toric variety defined by a 
reduction      fan $\pi^\vee (\Delta)_{red}$ and is of dimension less than $\dim X - \dim G$.
\end{rem}






\section{Virtual chambers and quotients}

Assume now that $X$ is projective and equipped with a $T$-linearized ample line bundle $L$.
This corresponds to a polytope $P$ in the dual lattice $M$ of $N$. Let $Q= \pi (P)$

\begin{cor}
Let $\widetilde{{\bf c}}$ be a virtual (real) cell of the projection $\pi:P \rightarrow Q$.
Then it defines a quotient toric variety
$X_{\widetilde{{\bf c}}}$  whose defining
fan is isomorphic to the fan in $(\ker \pi)^*$ induced by
the images of the locally (global) coherent costing $\Delta(\widetilde{{\bf c}})$. 
\end{cor}
\proof
This is a special case of Theorem \ref{mainquotientthm}.
\endproof

\section{Cox's quotient construction of toric varieties}

Let $\Delta$ be any fan in $N$ and $X_\Delta$ be the corresponding toric variety.
Cox has given a construction of $X_\Delta$ as the quotient of an open subset
of ${\Bbb C}^{|\Delta (1)|}$
where $\Delta (1)$ is the set of 1-dimensional cones of $\Delta$
(\cite{Cox95}). ${\Bbb C}^{|\Delta (1)}|$
comes equipped with a natural basis $\{ e_\rho \}$ where $\rho$ is the first lattice point
on an edge of $\Delta (1)$. 
There is canonical projection $\pi^\vee$ from  ${\Bbb C}^{|\Delta (1)|}$
to $N \otimes {\Bbb R}$ defined by sending each $e_\rho$ to $\rho$.
Let $G = \Hom_{\Bbb Z}(A_{n-1}(X_\Delta), {\Bbb C}^*)$.

\begin{cor} (Cox \cite{Cox95})
Let the notations be as above.
\begin{enumerate}
\item 
$\Delta$ canonically corresponds to 
a locally coherent costing of $\pi: {\Bbb C}^{|\Delta (1)|} \rightarrow N$
and the toric variety $X_\Delta$ is the quotient of ${\Bbb C}^{|\Delta(1)|}$ by the group $G$
defined by the above locally coherent costing.
\item $X_\Delta$ is a geometric quotient if and only if $\Delta$ is simplicial.
\end{enumerate}
\end{cor}
\proof 1. This is tautological.
Given $\sigma \in \Delta$, let $\sigma (1)$ be the set of edges of $\sigma$.
Clearly $\{ \gene_{\rho \in \sigma (1)} \{e_\rho\} | \sigma \in \Delta \}$ is the desired locally coherent string.

2. By Theorem \ref{mainquotientthm} (3), $X_\Delta$ is geometric if and only if $\Delta$ is tight.
It is tight if and only if the projection $\gene_{\rho \in \sigma (1)} \{e_\rho\} \rightarrow \gene_{\rho \in \sigma (1)} \{\rho\}$
does not drop dimension. This can happen if and only if $\gene_{\rho \in \sigma} \{\rho\}$
 is itself simplicial (i.e., $\{\rho\}_{\rho \in \sigma}$ are linearly independent).
\endproof

\section{Quotients defined by realizable virtual chambers}

Realizable virtual cells (i.e., deformed cells)
 are special. In this section, we shall explain that
their corresponding quotients share some characteristics of projective quotients defined by real cells. 
We will freely adopt the notations from \S \ref{deformcells}.

Let ${\bf c}_t$ and ${\bf c}_t'$ be two adjacent cells in $\Gamma_t$ such that the intersection
$${\bf c}_t \cap {\bf c}'_t={\bf c}^0_t$$ is a  face of each. Then we have 
$$U({\bf c}_t) \subset U({\bf c}^0_t)
\supset U({\bf c}'_t)$$ and the inclusions induce a diagram of morphisms of algebraic varieties
$$U({\bf c}_t)/\!/G  \;\;\;\;\;\;\;\;\;\;\; U({\bf c}'_t)/\!/G$$
$$f \searrow \;\;\;\;\;\;\;\;\;\;\; \swarrow f'$$
$$ U({\bf c}^0_t) /\!/G.$$

Set $\Sigma_0$ to be  $X^{{\rm ss}}({\bf c}^0_t)/\!/G \setminus  X^{\rm s}({\bf c}^0_t)/\!/G$.
Then $\Sigma_0$ admits a stratification by the so-called orbit types. Two points of $X$ have the same
orbit type if their stabilizers in $G$ are identical. This induces a $G$-invariant stratification of $X$
which in turn induces a stratification of $\Sigma_0$.

As in the projective case (e.g., $t=0$), we have

\begin{thm}
\label{1}
Assume that ${\bf c}^0_t$ is of codimension 1 and intersects with the interior of $\Gamma_t$. Then
\item{(i)}  $f$ and $f'$ are isomorphisms over the complement to $\Sigma_0$;
\item{(ii)} over each connected component $\Sigma_0'$ of a stratum of  $\Sigma_0 $, 
each fiber of  $f$ ($f'$) is isomorphic 
to a quotient of 
a weighted projective space of dimension $d$ ($d'$) by the finite group
 $\pi_0 (G_z)$ where $z$ is some  point in $X^{{\rm ss}}({\bf c}^0_t) \setminus  X^{\rm s}({\bf c}^0_t)$;
\item{(iii)} $d + d'  + 1 = {\rm codim} \; \Sigma_0'$.
\end{thm}
\proof  It is identical to the proof of Theorem 2.2 of \cite{Hu91}.
\endproof

\begin{rem}
\label{warning}
For  simplicity, we shall call the above birational transformation a Mori-type flips.
 Here we  firmly give a few warnings about the use of the term Mori-type flips.
\begin{enumerate}
\item Strictly, Mori flips are defined in the category of projective
varieties. We borrow the term in the non-projective case as well.
\item In a some special cases, some of our Mori-type flips may just correspond to blowups, 
while strictly Mori flips are birational transformations that are isomorphic in codimension 1.
\end{enumerate}
\end{rem}

\begin{prop}
Assume that ${\bf c}^0_t$ is of codimension 1,  lies on the boundary of $\Gamma_t$, and ${\bf c}_t$
is a chamber containing ${\bf c}^0_t$ as a face. Then
the map $f: X^{{\rm ss}}({\bf c}_t)/\!/G \rightarrow X^{{\rm ss}}({\bf c}^0_t)/\!/G$ is fiber bundle
whose typical fiber is isomorphic 
to a quotient of a weighted projective space by a finite abelian group.
\end{prop}

\proof
It is identical to the proof of Theorem 2.1 of \cite{Hu91} (the finite abelian group action was unfortunately
overlooked there).
\endproof

\begin{cor}
\label{rigid}
Every lexicographic (virtual) cell corresponds to a projective quotient which is,  
modulo  finite abelian group actions, 
a  tower of weighted projective bundles over a fixed point variety.
\end{cor}

It follows then

\begin{cor}
Every pair of quotient algebraic varieties
defined by realizable virtual cells are connected by a
 sequence of Mori-type flips.
\end{cor}

\proof
Assume that the two quotients are defined by a cell in $\Gamma_t$ and a cell $\Gamma_s$,
respectively. Both quotients can be flipped to 
the same projective quotient defined by
a lexicographic cell by crossing (virtual) walls. 
\endproof

\section{Some notes related to Chow quotient}

Besides the invariant theoretic quotients, there is a canonical ``quotient'' space, the Chow quotient.
To recall the construction (from \cite{KSZ}), note that the closure $\overline{G \cdot x}$ is a projective
subvariety, and as $x$ ranges within an $G$-invariant open subset $U \subset X$ of generic points,
these varieties will have the same dimension and degree. Let $\Ch(X)$ be the Chow variety of all 
algebraic cycles in $X$. The assignment $G \cdot x  \rightarrow \overline{G \cdot x}$
defines an embedding of $U/G$ into $\Ch(X)$. The closure of the image of $U/G$ in $\Ch(X)$
is the Chow quotient and denoted by $X/\!/^{ch} G$. The  Chow quotient is a toric
variety of the residue torus and its corresponding fan is the normal fan
of the fiber polytope $\Sigma (P,Q)$ (\cite{KSZ}).


Assume now that we are in the situation as in \S 13.
Then the $T$-linearized an ample line bundle over $X$ determines
a unique moment map $\mu_T$ for the $T$-action whose total
image is the polytope $P$.  It is a standard fact that for the image of  
$\overline{\cal O}_x = \overline{T \cdot x}$
 under $\mu_T$ is the face of $P$.
Also, with respect to the same ample line bundle $L$, a moment map of the subtorus $G$-action
is the composition $$\mu_G: X \rightarrow P \rightarrow Q$$
where the last map is the natural projection we began with.

Let $C= \sum_i a_i C_i  \in X/\!/^{ch}G$ where $C_i = \overline{G \cdot x_i}$ for some $x_i \in X$.
Then

\begin{prop}
$\{\mu_T ( \overline{\cal O}_{x_i})\}_i$ is a coherent string of $P \rightarrow Q$.
Moreover, every coherent string can be realized this way.
\end{prop}

\proof
This is a reformulation of Proposition 3.6 of \cite{KSZ}.
\endproof



It is known  that the Chow quotient maps to every projective GIT quotient.
This fact extends to non-projective cases as well.

\begin{cor}
Given any locally coherent costing $\Delta$,  there is a natural projection from the Chow quotient
$X/\!/^{ch}G$ to the quotient variety $U(\Delta)/\!/G$ that is a  birational $T/G$-toric morphism.
\end{cor}

\proof In the projective category, see the proofs in \cite{Kapranov} and \cite{KSZ}.

In our special situation, the following proof is quick: just observe from the definitions that
the fan of the Chow quotient refines the fan induced by any locally coherent costing.
The birationality and the $T/G$-equivariancy
is obvious. 
\endproof

\section{Bistellar flips and Mori-type flips}

 Let ${\cal A} = \{a_0, \ldots, a_n\}$ be a finite lattice point set in ${\Bbb Z}^d \subset {\Bbb R}^d$
and $Q=conv ({\cal A})$. We may consider $Q$ as the projection of the standard n-simplex $P$ by a suitable 
linear map $\pi$. In this case, the minimal elements of ($\cT_{coh}$) $\cT$ are the (coherent)
triangulations of $Q$ using solely the vertices from ${\cal A}$.

Coherent triangulations correspond to the vertices of the secondary polytope
$\Sigma({\cal A})$ of ${\cal A}$. Two vertices are joint by an edges of $\Sigma({\cal A})$ 
if and only if the corresponding triangulations are related by a bistellar flips. 
Here a bistellar flip is a local re-arrangement of a triangulation 
without introducing new vertices (see \cite{Ziegler}). 

 Embed ${\cal A} \subset 
{\Bbb R}^d \subset {\Bbb R}^{d+1}$ in the affine hyperplane of height
1 ($x_{d+1}=1$). Given any triangulation $T \in \cT$, taking the cone over $T$,
we get a fan $\Delta_T^0$ with the support $cone({\cal A})$. To get a complete fan $\Delta_T$, we add
the ray $R$ generated by the lattice point $0 \oplus (-1) \in 
{\Bbb R}^d \oplus {\Bbb R}$. The cones in $\Delta_T \setminus \Delta_T^0$ 
are generated by cones in $\partial
(cone({\cal A}))$ and the ray $R$.

\begin{prop}
$T$ is coherent if and only if $\Delta_T$ is projective.
\end{prop}

\proof 
The proof is almost tautological. Note that ($T$ is coherent) $\Delta_T$ is projective
if ($T$) $\Delta_T$ possesses a strictly convex piece-wise (affine) linear function.
Thus if $\Delta_T$ possesses such a function, restricting it to $T$ will give a desired
function for $T$. If $T$ has such a function, because $\Delta_T$ is obtained from
the cone over $T$ be adding one edge ${\Bbb R} (0 \oplus (-1))$, this function can be
extended to be a strictly convex piece-wise linear function
on $\Delta_T$ by assigning to $0 \oplus (-1)$ a generic number (to assure that
the extended linear functions are all different).
\endproof

\begin{exmp} Figure 4 is a non-coherent triangulation $T$ (\cite{Ziegler}) and 
its corresponding non-projective fan $\Delta_T$.
(cf. \cite{Fulton93} for another example of non-projective fans.) 
\end{exmp}

\vskip   1.8in
\begin{center}
Figure 4
\end{center}

Let  $\cT_{coh}^0$  be the subset of coherent triangulations that actually
use all vertices of ${\cA}$. 
Set
$$G_{\widehat{\cA}}=\Hom_{\Bbb Z}(A_{d}(X_{\Delta_T}), {\Bbb C}^*)$$
for any $T \in \cT_{coh}^0$. One checks that the group
does not depend on the choice of $T$. 
The group $ G_{\widehat{\cA}} $ acts on the complex vector space
$({\Bbb C})^{|\cA| +1}$ via the inclusion $G_{\widehat{\cA}} \subset ({\Bbb C}^*)^{|\cA |+1}.$

Now by our earlier results, we have

\begin{prop}
\begin{enumerate}
\item There is canonical correspondence between the triangulations
of $Q$ and quotient varieties of ${\Bbb C}^{|\cA|+1}$ by the group
$G_{\widehat{\cA}}$;
\item under this correspondence, coherent triangulations correspond
to projective quotient varieties, not-necessarily-coherent triangulations correspond
not-necessarily-projective quotient varieties;
\item  a bistellar flip between two triangulations transports to
a Mori-type flip\footnote{See Remark \ref{warning}.} between their corresponding quotients.
\end{enumerate}
\end{prop}

\section{Notes on virtual (pseudo)  oriented matroids}
This section is loosely related to  the rest of the paper,
especially \S\S 4--7. It raises more questions than it attempts to solve.

An oriented matroid of rank $n$ correspond to a (tamed) topological deformation  of an arrangement
of hyperplanes in ${\Bbb R}^n$ with that a realizable one corresponds to the hyperplane arrangement itself.
Given a hyperplane arrangement, it induces a (special) fan $\Delta$ (a normal fan of a zonotope) in ${\Bbb R}^n$.
For each hyperplane $e$, we assign an orientation: naming one side positive, the other negative.
Thus each cocell $\sigma$ (i.e., a cone in $\Delta$) corresponds to a sign vector $v^\sigma$ in $\{+,0,-\}^N$ according to
the location of the cocell relative to each and every hyperplane $e$,
where $N$ is the number of hyperplanes. So, $v^\sigma_e$ assumes one of the values $+, 0, -$
according to one of the three positions of the cocell $\sigma$ relative to the hyperplane $e$.

The oriented matroid defined by $\Delta$ is a combinatorial abstraction of properties of those sign vectors.
We refer the reader to $\cite{BLSWZ}$ for a list of axioms that characterizes the system of the sign vectors
$\{v^\sigma |\sigma \in \Delta\}$, which are  formal (and equivalent)
definitions of an oriented matroid. 
We use $\cM (\Delta)$ to denote this oriented matroid.
Such a $\cM (\Delta)$ is called a realizable oriented matroid.
 If we deform $\Delta$ to get $\Delta_t$,
the same strategy will still give us a sign vector for every  deformed cocell. Hence it also leads to
a (possibly non-realizable) oriented matroid $\cM (\Delta_t)$.

 In the above, we only consider the normal fan of a zonotope, a linear projection of a standard cube.
What about the normal fan of any polytope? Or more generally, what about an arbitrary fan?

For an arbitrary complete fan $\Delta$, it can always be extended to be the fan $\widehat{\Delta}$ induced by an arrangement of hyperplanes: just take the linear spans of the  cones of codimension one. This extension is canonical.

 So we still get a set of
hyperplanes and still can give orientations of these hyperplanes. Given a cone $\sigma$ in $\Delta$,
it may take a definite position relative to some hyperplane $e$: + side, $- $ side, or 0 (lying on the hyperplane).
For some other hyperplane $e$, $e$ runs though $\sigma$ in the relative interiors so that $\sigma$ does not
take any definite position relative to $e$. We may say that for $\sigma$, $e$ is a ``ghost'' hyperplane.
 Should this be problematic? What if we take this situation as a ``position''
too? So, let us use  the first letter in ``{\it u}ndecided'' to indicate this situation.
Then a cone $\sigma$ gives rise to a {\it generalized} sign vector 
$v^\sigma$ in $$\{+,0,-, u\}^N \setminus \{(u,\cdots,u)\}.$$
That is,  for any hyperplane $e$,
$v^\sigma_e$ assumes one of the values $+,0,-, u$ according to one of the above four positions.
Since $\Delta$ is not empty, for every given $\sigma$, $v^\sigma_e \ne u$ for some $e$, so that
$v^\sigma$ does not assume the value $\{(u,\cdots,u)\}.$ Now
we may call the combinatorial abstraction of the system of these generalized sign vectors a {\it virtual 
(or pseudo) oriented
matroid}, still denoted by $\cM (\Delta)$. Take any (tamed) topological  deformation $\Delta_t$ of $\Delta$,
we can assign a virtual oriented matroid $\cM (\Delta_t)$ to $\Delta_t$ in the same way as we did for $\Delta$ itself.
We may call the virtual oriented matroid $\cM (\Delta)$ a realizable one in accordance with the case
of oriented matroids.

Each generalized sign vector $v$ in $\{+,0,-, u\}^N$ splits into many (real) sign vectors:
whenever $v_e = u$  we can let it to split into to $\hat{v}_e = +, 0, -$. So, $v_e$ splits into
$3^{\chi (v)}$ many (real) sign vectors where $\chi (v)$ is the number of hyperplanes such that
$v_e =u$. The system of $\{\hat{v}_e\} $ is nothing but the system of the sign vectors for $\cM (\widehat{\Delta})$, where $\widehat{\Delta}$ is the canonical extension of $\Delta$ as explained
earlier.
We call $\cM (\widehat{\Delta})$ the oriented matroid extension of $\cM (\Delta)$. This extension is
canonical.

Hence given any system $\cM$ of {\it generalized} sign vectors, the above splitting procedure leads
to a canonical system $\widehat{\cM}$ of (real) sign vectors. Again, $\widehat{\cM}$ is called the canonical
extension of $\cM$. From this, we can give a practical definition of virtual oriented matroid as follows.

\begin{defn}
A set $\cM$ of {\it generalized} sign vectors
in $$\{+,0,-, u\}^N \setminus \{(u,\cdots,u)\}$$
 is called a {\it virtual } (or pseudo)  oriented matroid of rank $n$ if
its canonical extension  $\widehat{\cM}$ is an oriented matroid of rank $n$.
\end{defn}

An oriented matroid admits an involutive automorphism
$$I: v \to  - v.$$
A virtual oriented matroid does not have this property. In fact, if it does, it must be equal to
its canonical extension, hence itself an oriented
matroid. In this sense, 
{\it an oriented matroid is a virtual oriented matroid equipped with an involution}.
Put it slightly differently, a virtual oriented matroid is an ``oriented matroid without involution''.

The following diagram should be interesting,  which would extend the topological representation theorem of
oriented matroids to the virtual situation.
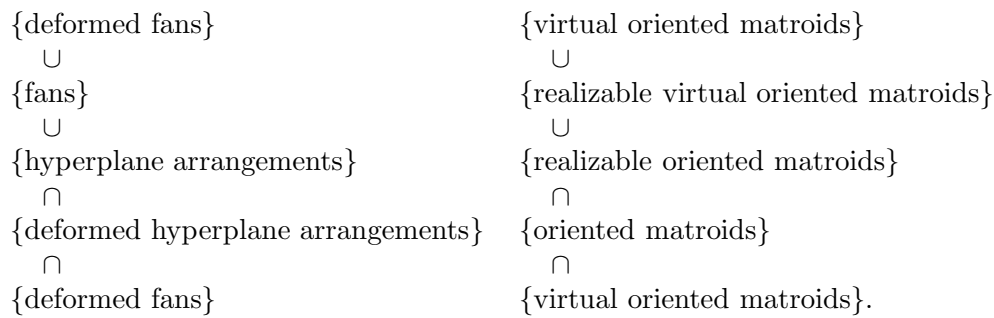
\begin{tabbing}  
 AAAAAAAAAAAAAAAAAAAAAAA  \= AAAAAAAAA \= haha \kill                      
\{\hbox{deformed fans}\} \>  \{\hbox{virtual oriented matroids}\} \\
$\;\;\;\;\cup $ \>  $\;\;\;\;\cup$\\
\{\hbox{fans}\}  \>   \{\hbox{realizable virtual oriented matroids}\} \\
$\;\;\;\;\cup$  \>  $\;\;\;\;\cup$\\
\{\hbox{hyperplane arrangements}\} \>   \{\hbox{realizable  oriented matroids}\} \\
$\;\;\;\;\cap $ \>  $\;\;\;\;\cap$\\
\{\hbox{deformed hyperplane arrangements}\}   \>     \{\hbox{oriented matroids}\}   \\
  $\;\;\;\;\cap $ \>  $\;\;\;\;\cap$\\
\{\hbox{deformed fans}\} \>  \{\hbox{virtual oriented matroids}\}. 
\end{tabbing}

We end our venture to virtual oriented matroids as follows.
 
\begin{conj}
With the help of the canonical oriented matroid extension, we may list a system of axioms which
characterizes a virtual oriented matroid in a way analogous to the case where a system of  axioms  characterizes
an oriented matroid. 
\end{conj}

\begin{prob}
Modify as many results of \cite{BLSWZ} as possible so that the revised results are valid
in the virtual situation. Moreover, find some interesting applications of the new theory.
\end{prob}

Finally, we would also like to ask

\begin{question}
Do virtual oriented matroids provide geometric structures for some combinatorial spaces
in analogous to that oriented matroids provide differentiable structures for MacPherson's
combinatorial differential manifolds \cite{MacP}?
\end{question}

\bibliographystyle{amsplain}
\makeatletter \renewcommand{\@biblabel}[1]{\hfill#1.}\makeatother


\end{document}